\documentclass[12pt]{article}

\usepackage{amsmath,amssymb,amsfonts,theorem,makeidx,latexsym,epsfig,subfigure}

\usepackage{color}

\newtheorem{defn}{Definition}[section]

\newtheorem{lemma}[defn]{Lemma}

{\theorembodyfont{\rmfamily}

\newtheorem{ex}[defn]{Example}}

\newtheorem{thm}[defn]{Theorem}

\newtheorem{prop}[defn]{Proposition}

\newtheorem{cor}[defn]{Corollary}

\numberwithin{equation}{section}

\newcommand{\ltr}{ L^2(\mathbb R) }

\newcommand{\mn}{\mathbb N}

\newcommand{\mr}{\mathbb R}

\newcommand{\mz}{\mathbb Z}

\newcommand{\mc}{\mathbb C}

\def\bp{{\noindent\bf Proof. \ }}

\def\ep{\hfill$\square$\par\bigskip}

\def\bqs{\begin{equation}}

\def\eqs{\tag*{$\square$}\end{equation}\par\bigskip}

\def\bop{\begin{op}\rm}

\def\eop{\end{op}}

\def\bee{\begin{eqnarray}}

\def\ene{\end{eqnarray}}

\def\bes{\begin{eqnarray*}}

\def\ens{\end{eqnarray*}}

\def\bei{\begin{itemize}}

\def\eni{\end{itemize}}

\def\bt{\begin{thm}}

\def\et{\end{thm}}

\def\bc{\begin{cor}}

\def\ec{\end{cor}}

\def\bpr{\begin{prop}}

\def\epr{\end{prop}}

\def\bl{\begin{lemma}}

\def\el{\end{lemma}}

\def\bd{\begin{defn}}

\def\ed{\end{defn}}

\def\bex{\begin{ex}}

\def\enx{\end{ex}}

\def\bfi{\begin{fig}}

\def\efi{\end{fig}}

\def\mzd{{\mathbb Z}^d}

\def\mrd{{\mathbb R}^d}

\def\bs{B^{\sharp}}

\newcommand{\ltrd}{ L^2({\mathbb R}^d) }

\def\bx{{\bf x}}
\def\by{{\bf y}}
\def\bn{{\bf n}}
\def\bk{{\bf k}}
\def\bj{{\bf j}}
\def\bm{{\bf m}}
\title{On partition of unities generated by entire functions
and Gabor frames in $\ltrd$ and $\ell^2(\mzd)$\thanks{This research was supported by Basic Science Research Program
through the National Research Foundation of Korea(NRF) funded by
the Ministry of Education(2013R1A1A2A10011922).
}}

\date{\today}

\author{Ole Christensen, Hong Oh Kim, Rae Young Kim}

\begin{document}

\maketitle

\begin{abstract} We characterize the entire functions $P$ of $d$ variables, $d\ge 2,$ for
which the $\mzd$-translates of $P\chi_{[0,N]^d}$ satisfy the partition of unity for
some $N\in \mn.$ In contrast to the one-dimensional case, these entire functions are not necessarily periodic.
In the case where $P$ is a trigonometric polynomial,
we characterize the maximal smoothness of $P\chi_{[0,N]^d},$ as well as the function that achieves it.
A number of especially attractive constructions are achieved, e.g., of trigonometric polynomials
leading to any desired (finite) regularity for a fixed support size.
As an application we obtain easy constructions of matrix-generated Gabor frames in $\ltrd,$ with
small support and high smoothness. By sampling this yields dual pairs of finite Gabor frames in $\ell^2(\mzd).$

\end{abstract}

\begin{minipage}{120mm}


\end{minipage}
\

\section{Introduction} \label{}
Partition of unity conditions appear in many different contexts in analysis, e.g.,
within harmonic analysis \cite{Fe4, FG2}. In this paper we characterize
the entire functions $P: \mc^d \rightarrow \mc, \, d\ge 2,$
which, for some fixed
$N\in\mn$, satisfy the partition of unity condition
\bee \label{157af-2N} \sum_{\bn \in \mzd} P(\bx + \bn) \chi_{[0,N]^d}(\bx + \bn)=1, \, \forall \bx \in \mrd.\ene
In contrast to the case $d=1$ treated in \cite{CKK6} such a function $P$ is not necessarily $(N\mz)^d$-periodic. In the special case of
periodic entire functions $P$ we derive an alternative and more direct characterization
of the partition of unity condition in terms of the Fourier coefficients of $P.$
For the case where $P$ is a trigonometric polynomial the  maximal smoothness of $P\chi_{[0,N]^d}$ is characterized,
as well as the entire functions $P$ that attain it.

The approach leads to a number of explicit constructions of
functions $P$ that yield a partition of unity and have desired
smoothness. We apply these results to provide very easy
constructions of pairs of dual Gabor frames with a number of
attractive properties, such as small support and high regularity.
Compared to the B-spline based frame constructions in \cite{CR}
these constructions are considerably more convenient: we avoid a
complicated ``book keeping," and the dual window has the same
support as the window itself. Due to the compact support and
continuity of the windows the results lead to an easy way to
construct finite dual Gabor frames in $\ell^2(\mzd)$ as well. For
more information on Gabor frames we refer to the monographs
\cite{G2,CB}.

The paper is organized as follows. In Section \ref{75} we characterize the entire
functions that have the partition of unity property, as well in the general case as in
the periodic case. The regularity issue is considered in Section \ref{41226a}, and the
applications to Gabor frames are in Sections \ref{41226b} and \ref{50401a}.

\section{Partition of unity for entire functions} \label{75}
Our first goal is to characterize
the entire functions $P: \mc^d \rightarrow \mc, \, d\ge 2,$
which, for some fixed
$N\in\mn$, satisfy the partition of unity condition \eqref{157af-2N}.
In order to do this we need to introduce some notation.
For $\by \in \mc^{d-1},$ write $\by=(y_1, \dots, y_{d-1}).$  For any $j\in \{1, \dots, d\},$
define the function
\begin{equation*}  
P_j: \mc \times \mc^{d-1}\to \mc, \, P_j(x, \by):=P(y_1, \dots, y_{j-1}, x, y_{j}, \dots, y_{d-1}).
\end{equation*}
The function $P_j$ is introduced in order to have a notation that allows us to
``pull out" the $j$th variable; in fact, for any $\bx=(x_1, \dots, x_d)\in \mc^d,$
letting $\by:= (x_1, \dots, x_{j-1}, x_{j+1}, \dots x_d)\in \mc^{d-1}$ yields that
\begin{equation*}
 P(\bx)= P(x_1, \dots, x_d)=P_j(x_j, (x_1, \dots, x_{j-1},
x_{j+1}, \dots x_d)) =P_j(x_j, \by).
\end{equation*}
 In particular,
$P(\bx)= P(x_1, \dots, x_d)=P_1(x_1,(x_2, \dots, x_d)).$
Given
$N\in \mn,$ the symbol \, $\sum_{\bn \in \mz_N^{d-1}}$ will denote
a sum over the $\bn =(n_1, \dots, n_{d-1})\in \mz^{d-1}$ for which all coordinates are
between $0$ and $N-1,$ i.e., \bes \sum_{\bn \in \mz_N^{d-1}}:=
\sum_{n_1=0}^{N-1} \dots\sum_{n_{d-1}=0}^{N-1}.\ens Now, for any
$j\in\{1, \dots, d\},$ define the function \bes Q_j: \mc \times
\mc^{d-1} \to \mc, \, Q_j(x, \by):=  \sum_{\bn \in \mz_N^{d-1}  }
P_j(x, \by + \bn).\ens Explicitly, \bes  Q_j(x, y_1, \dots,
y_{d-1})=   \sum_{n_1=0}^{N-1}
\dots\sum_{n_{d-1}=0}^{N-1} P(y_1+n_1, \dots, y_{j-1}+n_{j-1}, x,
y_{j}+n_{j}, \dots, y_{d-1}+n_{d-1}).\ens
We will now characterize the entire functions $P$ satisfying the partition of unity condition in terms of (any of) the associated functions $Q_j.$
\bpr \label{157d-2N} Let $P: \mc^d \to \mc$ be an entire function. Let $N\in \mn,$ and consider any $j\in \{1, \dots, d\}.$
Then the following are equivalent:
\begin{itemize}
    \item[{\rm (a)}] $P$  satisfies the partition of unity condition \eqref{157af-2N};
    \item[{\rm (b)}] For any $\by\in[0,1]^{d-1},$
    the restriction of $Q_j(\cdot, \by)$ to $\mr$ is
    $N$-periodic and the Fourier coefficients $c_k(\by)$ in the expansion
\bee \label{157c-2N}  Q_j(x,\by)=\sum_{k\in \mz} c_k(\by) e^{2\pi ikx/N}, \ x\in \mr  \ene
    satisfy that  $c_k(\by) = \frac1{N} \delta_{k,0}$ for $k\in N\mz.$
\end{itemize}
\epr
\bp (a)$\Rightarrow$(b)
Assume first that \eqref{157af-2N} holds, and take any $j\in \{1, \dots, d\}.$ Then
$\sum_{\bn \in \mz_N^d}
P(\bx + \bn) =1, \, \forall \bx\in [0,1]^d.$
Thus, for any fixed  $\by \in [0,1]^{d-1},$ \bee \label{108a-2N}
\sum_{\ell=0}^{N-1} Q_j(x+\ell,\by)=
\sum_{\ell=0}^{N-1}\sum_{\bn \in \mz_N^{d-1}  } P_j(x+\ell, \by +
\bn)=1, \, x\in [0,1].\ene Since $Q_j(\cdot,\by)$  is an entire
function, \eqref{108a-2N} then holds for all $x \in \mr.$
Replacing $x$ by $x+1$ in \eqref{108a-2N} and
 subtracting the two expressions shows that
$Q_j(x+N,\by)=Q_j(x,\by), \  x\in \mr.$
So we conclude that the restriction of
$Q_j(\cdot,\by)$ to $\mr$ is $N$-periodic. Writing $Q_j(\cdot,\by)$ as the Fourier series \eqref{157c-2N},
the equation \eqref{108a-2N} takes the form
\begin{equation} \label{et-65-2N}
 \sum_{k\in \mz} c_k(\by)\,\left[1+ e^{2\pi i k/N}
+ \cdots + \left( e^{2\pi i k/N}\right)^{N-1} \right] \, e^{2\pi ikx/N}=1.
\end{equation}
We note that
\begin{equation} \label{et-66-2N}
1+ e^{2\pi i k/N}
+ \cdots + \left( e^{2\pi i k/N}\right)^{N-1}=\left\{
\begin{array}{ll}
   N, & k\in N\mz \\
   0, & k\notin N\mz.
\end{array}
\right.
\end{equation}
From \eqref{et-65-2N} and \eqref{et-66-2N}, we see that
$c_k(\by)=\frac{1}{N}\delta_{k,0}$ for $k\in N\mz,$  as claimed.

\noindent (a)$\Leftarrow$(b) Let again $j\in \{1, \dots, d\},$ and
consider any $\bx\in [0,1]^d.$ Then
\bes & \ &  \sum_{\bn\in \mz^d} P(\bx + \bn) \chi_{[0,N]^d}(\bx + \bn)
=   \sum_{\bn\in \mz_N^d} P(\bx + \bn)\\
&=&
  \sum_{n_1=0}^{N-1}  \sum_{n_2=0}^{N-1} \cdots \sum_{n_d=0}^{N-1}
P_j(x_j+n_j, (x_1+n_1, \dots, x_{j-1}+n_{j-1}, x_{j+1}+n_{j+1},
\dots, x_d+n_d)). \ens
With
$\by:=(x_1+n_1, \dots, x_{j-1}+n_{j-1}, x_{j+1}+n_{j+1}, \dots,
x_d+n_d),$  the assumption in (b) and
\eqref{et-66-2N} now yields that \bes  & \ &  \sum_{\bn\in \mz^d}
P(\bx + \bn) \chi_{[0,N]^d}(\bx + \bn) =   \sum_{n_j=0}^{N-1}
Q_j(x_j+n_j, \by)  = \sum_{n_j=0}^{N-1}\sum_{k\in \mz}
c_k(\by) e^{2\pi ik(x_j+n_j)/N} \\ & = & \sum_{k\in \mz}
c_k(\by)\,\left[1+ e^{2\pi i k/N} + \cdots + \left( e^{2\pi i
k/N}\right)^{N-1} \right] \, e^{2\pi ikx_j/N} = 1. \ens By periodicity of $\sum_{\bn\in \mz^d} P(\cdot +\bn) \chi_{[0,N]^d}(\cdot + \bn)$,
  \eqref{157af-2N} therefore
holds for all $\bx\in \mrd.$ \ep

Note that if the conditions in Proposition \ref{157d-2N} hold, then (b) actually holds for all $\by \in \mrd.$
From the proof we also see immediately
that a similar result holds with the ``square" $[0,N]^d$ replaced by a
rectangle  $[0, N_1] \times \dots \times [0, N_d],$ where $N_1, \dots, N_d\in \mn.$

In  \cite{CKK6} it was proved that in the case $d=1,$ an entire function
$P:\mc \to \mc$ satisfying \eqref{157af-2N} is automatically  $N$-periodic. The following example shows that this does not generalize to the case $d>1.$
In fact, an entire function $P: \mc^d \to \mc$
satisfying \eqref{157af-2N}  might not be periodic in any of the variables:

\begin{ex} \label{76} Let $d=2$ and let
$f: \mc \to \mc$ denote an entire function.
Consider the entire function $P: \mc^2\to \mc$ given by \bes P(x_1,x_2):=e^{\pi i x_2}
f(x_1) +\frac{1}{4}+e^{\pi i x_1}f(x_2).\ens An easy direct computation shows that
$P$ satisfies the partition of unity condition \eqref{157af-2N} for $N=2.$
Alternatively, let $j=1$ and fix $x_2\in [0,1]$. Then
\begin{eqnarray*}
   Q_1(x_1,x_2)= P(x_1,x_2)+P(x_1,x_2+1)
           = \frac12 + e^{\pi i x_1}\left(f(x_2)+f(x_2+1)\right).
\end{eqnarray*}
Thus, for any fixed $x_2,$ $Q_1(\cdot, x_2)$ is 2-periodic  and
satisfies the condition (b) in
Proposition \ref{157d-2N}; this again
implies  that $P$ satisfies the partition of unity condition
\eqref{157af-2N}. However, in generel  $P(\cdot, \cdot)$ is
periodic neither
in the first variable nor in the second variable.\ep
\end{ex}

In order to have an extra technical tool at our disposal (namely, Fourier series)  we
will now restrict our attention to entire functions $P: \mc^d \to \mc$ that
are  $(N\mz)^d$-periodic, i.e.,
entire functions $P$ for which the restriction to $\mrd$ can be written  in the form
\begin{equation} \label{77}
P(\bx):=\sum_{\bk\in\mzd} c_{\bk} e^{2\pi i \bk \cdot \bx /N}, \, \bx\in \mrd.
\end{equation}
For this class of entire functions we will now give a convenient characterization of
the functions having the partition of unity property in terms of the Fourier coefficients $c_\bk.$
The reader who checks the proof will notice that a
similar result holds for entire functions that are periodic along a lattice $N_1\mz \times \cdots \times N_d \mz$ for some $N_1, \dots, N_d\in \mn.$

\bc \label{eh-4} An entire $(N\mz)^d$-periodic function $P$
of the form \eqref{77} satisfies \eqref{157af-2N} if and only if
\begin{equation}\label{50104k}
   c_{\bk}=\frac{1}{N^d}\delta_{\bk,{\bf 0}}, \, \forall \bk\in (N\mz)^d,
\end{equation}
\ec

\bp We will apply
Proposition \ref{157d-2N} for the choice $j=1.$
For $x\in \mr$  and $\by \in [0,1]^{d-1},$ we have
\bee \label{79}
Q_1(x,\by)&=&\sum_{\bn\in \mz_N^{d-1}} P_1(x,\by+\bn)
 =  \sum_{\bn\in \mz_N^{d-1}} P(x,\by+\bn). \ene
Let us write the Fourier series \eqref{77} in a slightly  different form, namely, as
\begin{equation*} 
P(\bx)=P(x_1, \tilde{\bx})=\sum_{k\in\mz} \sum_{\bm \in \mz^{d-1}} c_{k, \bm} e^{2\pi i kx_1 /N}
 e^{2\pi i \bm \cdot \tilde{\bx}/N},
\end{equation*} where $\tilde{\bx}:=(x_2, \dots, x_d).$ Inserting this in \eqref{79} yields
\bes
Q_1(x,\by)&=& \sum_{\bn\in \mz_N^{d-1}}\sum_{k\in\mz} \sum_{\bm \in \mz^{d-1}} c_{k, \bm} e^{2\pi i kx /N}
 e^{2\pi i \bm \cdot (\by + \bn)/N}  \\ & = &
\sum_{k\in\mz}   \left( \sum_{\bn\in \mz_N^{d-1}} \sum_{\bm \in \mz^{d-1}} c_{k, \bm}
 e^{2\pi i \bm \cdot (\by + \bn)/N}\right) e^{2\pi i kx /N}   =
\sum_{k\in \mz} c_k(\by) e^{2\pi i kx /N},\ens
where the coefficients $c_k(\by)$ are given by
\bee \notag  c_k(\by) & = & \sum_{\bn\in \mz_N^{d-1}} \sum_{\bm \in \mz^{d-1}} c_{k, \bm}
 e^{2\pi i \bm \cdot (\by + \bn)/N}  = \sum_{\bm \in \mz^{d-1}} c_{k, \bm}
 e^{2\pi i \bm \cdot \by/N} \sum_{\bn\in \mz_N^{d-1}} e^{2\pi i \bm \cdot \bn/N}.\ene
 The sum
$\sum_{\bn\in \mz_N^{d-1}} e^{2\pi i \bm \cdot \bn/N}$
is only nonzero whenever $\bm = N{\bf p} $ for some ${\bf p} \in \mz^{d-1},$ in which case the sum is $N^{d-1};$
inserting this yields
\bee \label{81} c_k (\by) & = & N^{d-1} \sum_{{\bf p} \in \mz^{d-1}} c_{k, N{\bf p}}
 e^{2\pi i {\bf p} \cdot \by}. \ene

The condition (b) in  Proposition \ref{157d-2N}, namely, that the coefficients $c_k(\by)$
in \eqref{81} satisfy that
$c_k(\by)=\frac{1}{N}\delta_{k,0}$ for $k\in N\mz,$ is equivalent to
$$c_{0, N{\bf p}}=\frac1{N^d}\delta_{{\bf 0}, {\bf p}}, \, \mbox{ for } {\bf p} \in \mz^{d-1};
\ c_{k, N{\bf p}}=0, \, \, \mbox{ for } k\in N\mz\setminus\{0\},
{\bf p} \in \mz^{d-1},
$$
which again is equivalent with $ c_{\bk}=\frac{1}{N^d}\delta_{\bk,0}, \, \forall \bk\in (N\mz)^d.$
This completes the proof.
\ep

Note that by Corollary  \ref{eh-4} the entire $(N\mz)^d$-periodic functions $P$
in the form \eqref{77} which satisfies \eqref{157af-2N}, precisely are the ones having the form
$P(\bx)=\frac{1}{N^d} + \sum_{\bk\in \mzd\setminus (N\mz)^d} c_{\bk}e^{2\pi i \bk \cdot\bx/N}.$
In general such functions are not tensor products of $d$ one-dimensional Fourier series.

Our next goal is to construct entire functions $P$ such that $P\chi_{[0,N]^d}$ satisfies the partition of unity condition and has desired
regularity. Clearly, already the continuity of $P\chi_{[0,N]^d}$ forces $P$ to vanish
along the boundary of the ``square"  $[0,N]^d.$ The following example shows that even this extra constraint does not imply that $P$ is periodic:

\bex \label{176} Let $d=2$ and $N=2.$ Our purpose is to construct a nonperiodic entire function which satisfies the
partition of unity condition and vanishes on the
boundary of $[0,2]^2.$  Let $f,g$ be any entire functions. A calculation like in Example \ref{76} shows that the  entire function
\bes P(x_1,x_2):=\frac{1}{4}+e^{\pi i x_2}
f(x_1) +e^{\pi i x_1}f(x_2)+ e^{-\pi i x_2}
g(x_1) +e^{-\pi i x_1}g(x_2)  \ens
satisfies the partition of unity condition.
Now choose $f$ as an entire {\it nonperiodic} function such that $f(0)=f(2)=0,$
and let $g(x):= -\frac14 -f(x)+ \frac18 e^{-\pi ix}.$ Then
$g(0)=g(2)=-1/8.$ Also, for any $x_1, x_2\in \mr,$
$P(x_1,0) =
P(x_1,2)=
P(0, x_2)=
P(2, x_2) =0.$
Thus, the function $P$ vanishes on the boundary of $[0,2]^2.$ Note that
\bes P(x_1,x_2) & = & \frac14  +e^{\pi i x_2}
f(x_1) +e^{\pi i x_1}f(x_2) + e^{-\pi i x_2}
\left(  -\frac14 -f(x_1)+ \frac18 e^{-\pi ix_1} \right) \\ & \ & +e^{-\pi i x_1} \left(  -\frac14 -f(x_2)+ \frac18 e^{-\pi ix_2} \right) \\ & = & \frac14 + f(x_1) \left[ e^{\pi ix_2}-e^{-\pi ix_2}   \right]+ f(x_2) \left[ e^{\pi ix_1}-e^{-\pi ix_1}   \right]   \\ & \ &   -\frac14 e^{-\pi ix_1}  -\frac14 e^{-\pi ix_2} +  \frac14 e^{-\pi ix_1}  e^{-\pi ix_2},\ens which, by the assumptions on $f,$ clearly is nonperiodic.
\ep \enx

The Fourier series turn out to be the key ingredient in the regularity discussion in the next section, so we will continue to assume that $P$ is periodic.

\section{Regularity } \label{41226a}

Among the $(N\mz)^d$-periodic entire functions $P: \mc^d\to \mc$
satisfying the partition of unity condition, we will now restrict
our attention to the ones with only a finite number of nonzero
Fourier coefficients, i.e., the trigonometric polynomials. For
such functions, the following result characterizes the maximal
smoothness of $P\chi_{[0,N]^d},$ as well as the entire functions
$P$ that achieve it. Recall that for $L\in \mn,$ the space
$C^{L-1}(\mrd)$ consists of the functions $f: \mrd\to \mc$ for which all
the partial derivatives
$$ \frac{\partial^\ell f}
{\partial x_1^{\ell_1}\partial x_2^{\ell_2}\cdots\partial x_d^{\ell_d}}, \, \, 0\leq \ell=\ell_1+\cdots+\ell_d \leq L-1,
$$
exist and are continuous.

\begin{thm}\label{eh-3} Let $K,N\in\mn$. Assume that
\begin{equation} \label{77-N}
P(\bx)=\sum_{\bk \in \left(\mz \cap[-K,K]\right)^d} c_{\bk}
e^{2\pi i \bk \cdot \bx /N},\ \bx\in\mr^d
\end{equation}
is a real-valued trigonometric polynomial.
 Then the following hold:
\begin{itemize}
 \item[{\rm (a)}]  There does not exist $P$ of the form \eqref{77-N} such that $P\chi_{[0,N]^d}\in
    C^{2K}(\mrd)$;
    \item[{\rm (b)}] Fix $L \in \{1,2,\cdots, 2K \}.$ Then  $P\chi_{[0,N]^d} \in C^{L-1}(\mr^d)$ if and only if
\begin{equation}\label{eh-18}
P(\bx)=P(x_1,\cdots, x_d)= \prod_{j=1}^d \left(e^{\pi i x_j/N} \sin(\pi
x_j/N)\right)^{L} A_{L}(\bx)
\end{equation}
 for a  trigonometric polynomial
\begin{equation}\label{et-69}
   A_{L}(\bx)= \sum_{\bk \in \left(\mz \cap[-K,K-L]\right)^d} a_{\bk} e^{2\pi i \bk\cdot \bx/N};
\end{equation}
 \item[{\rm (c)}] Assume that $P$ has the form \eqref{eh-18}. Then $P\chi_{[0,N]^d}$  satisfies the partition of unity
condition \eqref{157af-2N} if and only if for all $\bm\in (N\mz)^d,$
\begin{equation}\label{eh-17}
   \sum_{\scriptsize
\begin{array}{c}
   \bj+\bk=\bm \\
   \bj\in \left(\mz \cap[0,L]\right)^d \\
   \bk \in \left(\mz \cap[-K,K-L]\right)^d
\end{array}}
(-1)^{Ld-(j_1+\cdots+ j_d)}
   {L\choose j_1}\cdots {L\choose j_d}a_{\bk}
   =
   \left(\frac{(2i)^L}{N}\right)^d \delta_{\bm,{\bf 0}}.
\end{equation}
\end{itemize}
\end{thm}
\bp We will first prove (b). First  note that
if $P\chi_{[0,N]^d}\in
C^{L-1}(\mrd)$, then  for $j=1,\dots,d$, we have
$\frac{\partial^\ell P}{\partial x_j^\ell}(x_1,\cdots, x_{j-1},0,x_{j+1},\cdots,x_d)=0,\ \ell=0,\cdots,L-1.$
Let us write \eqref{77-N} as
$P(\bx)=P(x_1,\cdots,x_d)=\sum_{k_1=-K}^{K} c_{k_1}(x_2,\cdots,x_d)e^{2\pi i k_1 x_1/N},$
where
$$ c_{k_1}(x_2,\cdots,x_d)=\sum_{k_2=-K}^{K} \cdots \sum_{k_d=-K}^{K}
c_{\bk}e^{2\pi i (k_2,\cdots,k_d)\cdot (x_2,\cdots,x_d) /N}.$$
By \cite[Theorem 3.1]{CKK6},
we know that whenever $L\le 2K,$ then $\frac{\partial^\ell P}{\partial x_1^\ell}(0,x_2,\cdots,x_d)=0$ for $ \ell=0,\dots,L-1$ if and only if
\begin{equation}\label{eh-6}
P(\bx)= \left(e^{\pi i x_1/N} \sin(\pi
x_1/N)\right)^{L} A_{1,L}(\bx)
\end{equation} for a  trigonometric polynomial
\begin{equation} \label{eh-7}
   A_{1,L}(\bx):= \sum_{k_1=-K}^{K-L} a_{k_1}(x_2, \cdots, x_d) e^{2\pi i k_1x_1/N}.
\end{equation}
We note that for some coefficients $d_{\bk},$
\begin{eqnarray}
a_{k_1}(x_2, \cdots, x_d)
&=&\sum_{k_2=-K}^{K} \cdots \sum_{k_d=-K}^{K} d_{\bk} e^{2\pi i (k_2,\cdots, k_d)\cdot (x_2,\cdots,x_d)} \nonumber \\
&=&\sum_{k_2=-K}^{K}\left(\sum_{k_3=-K}^{K}\cdots
\sum_{k_d=-K}^{K} d_{\bk}e^{2\pi i (k_3,\cdots, k_d)\cdot
(x_3,\cdots,x_d)}\right)e^{2\pi i k_2 x_2}. \label{eh-30}\end{eqnarray}
Now, by
\eqref{eh-6} and \eqref{eh-7}, $\frac{\partial^\ell P}{\partial
x_2^\ell}(x_1,0,x_3,\cdots,x_d)=0$ for $0\leq \ell\leq L-1$ if and
only if $\frac{\partial^\ell A_{1,L}}{\partial
x_2^\ell}(x_1,0,x_3,\cdots,x_d)=0$ for $0\leq \ell\leq L-1;$ or, if
and only if $\frac{\partial^\ell a_{k_1}}{\partial
x_2^\ell}(0,x_3,\cdots,x_d)=0$ for $0\leq \ell\leq L-1.$
By \cite[Theorem 3.1]{CKK6} and \eqref{eh-30} again, this is equivalent to
$$a_{k_1}(x_2,\cdots,x_d)= \left(e^{\pi i x_2/N} \sin(\pi
x_2/N)\right)^{L}\sum_{k_2=-K}^{K-L} a_{k_1,k_2}(x_3,\cdots, x_d) e^{2\pi i k_2x_2/N}.$$
That is,
$P(\bx)= \prod_{j=1}^2\left(e^{\pi i x_j/N} \sin(\pi
x_j/N)\right)^{L} A_{2,L}(\bx)$
for a trigonometric polynomial
$$
A_{2,L}(\bx)=\sum_{k_1=-K}^{K-L}
\sum_{k_2=-K}^{K-L} a_{k_1,k_2}(x_3,\cdots, x_d) e^{2\pi i (k_1x_1+k_2x_2)/N}.$$
Inductively,
we have
$P(\bx)= \prod_{j=1}^d \left(e^{\pi i x_j/N} \sin(\pi
x_j/N)\right)^{L} A_{L}(\bx)$ for a  trigonometric polynomial
\begin{equation*}
   A_{L}(\bx)= \sum_{k_1=-K}^{K-L}\cdots \sum_{k_d=-K}^{K-L} a_{k_1,\cdots,k_d} e^{2\pi i (k_1,\cdots,k_d) \cdot (x_1,\cdots, x_d)/N}.
\end{equation*}

Conversely, assume that
$P(\bx)= \prod_{j=1}^d \left(e^{\pi i x_j/N} \sin(\pi
x_j/N)\right)^{L} A_{L}(\bx)$ for a  trigonometric polynomial $A_L$
of the form \eqref{et-69}.
Then
for $0\leq \ell_1+\cdots+\ell_d \leq L-1$, and for $1\leq m\leq d$,
$$ \frac{\partial^\ell P}
{\partial x_1^{\ell_1}\partial x_2^{\ell_2}\cdots\partial x_d^{\ell_d}}
(x_1, \cdots,x_{m-1}, 0,x_{m+1},\cdots, x_d)=0
$$
is trivially satisfied.
Since $P$ is entire and
$(N\mz)^d$-periodic, $P\chi_{[0,N]^d} \in C^{L-1}(\mr^d)$. Hence
 (b) holds.

In order to prove (a), assume that $P\chi_{[0,N]} \in C^{2K}(\mr^d)$.
Using (b) with $L=2K,$ we see that $P(\bx)=\alpha \prod_{j=1}^d \left( \sin(\pi
x_j/N)\right)^{2K}$ for some $\alpha\in \mc.$
A direct calculation shows that for $(x_2,\dots x_d)\in [0,N]^{d-1}$,
$$\dfrac{\partial^{2K}P}{\partial x_1^{2K}}(0,x_2\cdots x_d)
=\alpha \left(\frac{\pi}{N}\right)^{2K}(2K)! \prod_{j=2}^d \left( \sin(\pi
x_j/N)\right)^{2K}
\not\equiv 0. $$
This is a contradiction. Thus $P\chi_{[0,N]^d} \not\in C^{2K}(\mr^d)$ and (a) holds.

For the proof of (c), we use the identity
\begin{eqnarray*}
\left(e^{\pi i x/N} \sin(\pi x/N)\right)^{L} =
\left( \frac{e^{2\pi i x/N} -1}{2i} \right)^L
= \left( \frac{1}{2i} \right)^L
\sum_{j=0}^L { L \choose j} e^{2\pi i j x  /N} (-1)^{L-j};
\end{eqnarray*}
then
\begin{eqnarray*}
  P(\bx)&=&\left( \frac{1}{2i} \right)^{Ld}
\sum_{\bj\in \left(\mz \cap[0,L]\right)^d} { L \choose j_1}\cdots{ L \choose j_d}
e^{2\pi i \bj \cdot \bx /N} (-1)^{Ld-(j_1+\cdots+j_d)}
\sum_{   \bk \in \left(\mz \cap[-K,K-L]\right)^d
} a_{\bk} e^{2\pi i \bk\cdot \bx/N}\\
&=&\left( \frac{1}{2i} \right)^{Ld} \sum_{\bm\in \mzd}
\left(\sum_{
\scriptsize
\begin{array}{c}
   \bj+\bk=\bm \\
   \bj\in \left(\mz \cap[0,L]\right)^d \\
   \bk \in \left(\mz \cap[-K,K-L]\right)^d
\end{array}
}
   (-1)^{Ld-(j_1+\cdots+ j_d)}
   {L\choose j_1}\cdots {L\choose j_d}a_{\bk}\right)e^{2\pi i \bm\cdot \bx/N}.
\end{eqnarray*}
By Corollary \ref{eh-4}, the condition that
$P\chi_{[0,N]^d}$  satisfies the partition of unity
condition is equivalent to \eqref{eh-17}.
Hence (c) holds. \ep

Let us use Theorem \ref{eh-3} to construct a partition of unity explicitly.

\begin{ex}
   Let $L=N=K=d=2.$ We will find $P$ of the form \eqref{eh-18}
such that $P\chi_{[0,2]^2} \in C^1(\mr^2)$ and
the partition of unity condition \eqref{157af-2N} holds.
Let $A_2(x_1,x_2)= \sum_{k_1, k_2 \in \{-2,-1,0\}}
a_{k_1,k_2}  e^{\pi i (k_1x_1 +k_2 x_2)},$ and assume that the coefficients $a_{k_1,k_2}$
are real numbers and satisfy that
\begin{equation*} 
a_{k_1,k_2}=a_{-k_1-2,-k_2-2}, \ \ k_1,k_2=-2,-1,0.
\end{equation*}
We now apply Theorem \ref{eh-3} (c).
Due to  the limitations on the summation indices ${\bf j,k}$
in \eqref{eh-17}, it is now enough to choose the coefficients $a_{k_1,k_2}$ such that \eqref{eh-17} is satisfied
for $\bm=(0,0),\pm(2,0),\pm(0,2),\pm(2,2).$
A direct calculation  shows that
these equations amount to the four equations $2 a_{0,0}-4 a_{-1,0} +2 a_{0,-2} -4 a_{0,-1} +4 a_{-1,-1} =4 ; \ \ a_{0,-2} -2 a_{0,-1} +a_{0,0}=0 ; \ \
a_{0,-2} -2 a_{-1,0} +a_{0,0}=0 ; \  \
a_{0,0} =0.$
If we assign a parameter $a_{0,-2}=t, \, t\in \mr,$ the solution can be written in vector form as
\begin{equation*}
\begin{pmatrix}
a_{0,0} \\
a_{-1,0} \\
a_{-1,-1} \\
a_{0,-1} \\
a_{0,-2} \
\end{pmatrix}
=\begin{pmatrix}
a_{-2,-2} \\
a_{-1,-2} \\
a_{-1,-1} \\
a_{-2,-1} \\
a_{-2,0} \
\end{pmatrix}
=t
\begin{pmatrix}
0 \\
1/2\\
1/2\\
1/2\\
1
\end{pmatrix}
+\begin{pmatrix}
0 \\
0\\
1\\
0\\
0
\end{pmatrix}.
\end{equation*}
By direct calculation using  \eqref{eh-18} it now follows that for any $t\in \mr,$ the
trigonometric polynomial
\begin{eqnarray*}
  P(x_1,x_2) & = & \sin^2(\pi x_1/N) \sin^2(\pi x_2/N)
\left[t\left(
\cos(\pi x_1) +  \cos(\pi x_2)  + 2  \cos(\pi (x_1-x_2)) +1/2\right)
+1\right]
\end{eqnarray*}
satisfies the requirements.
\ep

\end{ex}

Let us comment on the parameters $K,N,L$ appearing in Theorem \ref{eh-3}.
Not surprisingly, Theorem \ref{eh-3} shows that the ``budget of
nonzero Fourier coefficients" for the entire function $P$ in \eqref{77-N} limits
the possible smoothness of $P\chi_{[0,N]^d}.$
If we
fix $N\in \mn,$ e.g., if we want a certain support size for the function $P\chi_{[0,N]^d},$ the condition \eqref{50104k} is automatically satisfied for
$\bk \neq {\bf 0}$ if we take $K\le N-1$ in \eqref{77-N}. Thus, the characterization in Theorem \ref{eh-3} (b) yields an easy way to construct partition of unities with smoothness
at most $2(N-1)-1= 2N-3.$ We  note that if we for some $N\in \mn$
want maximal smoothness of $P\chi_{[0,N]^d},$ there is a unique function $P$
among the functions in \eqref{77-N} with $K\le N-1$ that achieves it. In fact,
fixing $N\in \mn$ and taking $K=N-1, L=2K,$ yields that
$K-L=-K.$ Thus $P$ must have the form
\bes \notag P(\bx) & = & \prod_{j=1}^d \left(e^{\pi i x_j/N} \sin(\pi
x_j/N)\right)^{2K} \sum_{\bk \in \left(\mz \cap[-K,-K]\right)^d} a_{\bk} e^{2\pi i \bk\cdot \bx/N};\ens  that is, for some $a\in \mr,$
\bee  P(\bx) & = & a\, \prod_{j=1}^d e^{2\pi i Kx_j/N} \sin^{2K}(\pi
x_j/N)   e^{-2\pi i(Kx_1+\cdots +Kx_d)/N} = a\, \prod_{j=1}^d \sin^{2N-2}(\pi x_j/N). \, \, \, \, \,
\label{50104q} \hspace{.2cm} \ene

This construction stands out as the optimal one with regard to smoothness, and also as the  simplest and most elegant one. Let us formulate the result formally:

\bc \label{eh-5}  Let $N\in\mn$ and
let \bee \label{41228g} P(\bx):=    \left( \frac{4^{N-1}}{N {2N-2
\choose N-1}}\right)^d   \prod_{j=1}^d \sin^{2N-2}(\pi x_j/N). \ene
Then $P\chi_{[0,N]^d}$  satisfies the partition of unity
condition \eqref{157af-2N} and belongs to $C^{2N-3}(\mr^d)$. \ec
\bp The result follows almost immediately from \eqref{50104q}; we just have to determine
the value of $a$ such that \eqref{50104k} is  satisfied for
$\bk ={\bf 0}.$ Using that
\begin{equation*}
  \sin^{2N-2}(\pi x/N)=
  \left( \frac{e^{\pi  i x /N} -e^{-\pi i x/N}}{2i} \right)^{2N-2}
 = \frac{1}{4^{N-1}} \sum_{k=-N+1}^{N-1} (-1)^{k} {2N-2\choose N-1+k}
 e^{2\pi i k x/N}
 \end{equation*}
and pulling out the coefficient corresponding to $k=0,$ leads to the form in
\eqref{41228g}.
\ep

In Corollary \ref{eh-5} the regularity of $P\chi_{[0,N]^d}$ is related to the support size, i.e., to the parameter $N.$ At the price of a more complicated construction, arbitrary (finite) regularity can be obtained even for $N=2.$

\bc \label{eh-15-3} Given $L\in \mn,$  let
\bes P(\bx)= \prod_{j=1}^d  \sin^{2L}(\pi x_j/2) \sum_{k_j=0}^{L-1}{2L-1 \choose k_j}
   \sin^{2(L-1 -k_j)}(\pi x_j/2) \cos^{2k_j} (\pi x_j/2).\ens
   Then $P\chi_{[0,2]^d}$ satisfies the partition of unity condition and belongs to
   $C^{2L-1}(\mrd).$  \ec

\bp In the case $d=1,$ it was shown in \cite{CKK6} that letting
\begin{equation*}
   Q(x):=\sin^{2L}(\pi x/2) \sum_{k=0}^{L-1}{2L-1 \choose k}
   \sin^{2(L-1 -k)}(\pi x/2) \cos^{2k} (\pi x/2), \, x\in \mr,
\end{equation*} we have that
$Q\chi_{[0,2]}$ belongs to $C^{2L -1}(\mr)$ and satisfies the partition of unity property.
Since $P(\bx)= \prod_{j=1}^d Q(x_j),$ we have
\begin{eqnarray*}
\sum_{\bn\in \mz_N^d} P(\bx + \bn)&=&
\sum_{n_1=0}^{N-1}  \cdots \sum_{n_d=0}^{N-1}
Q(x_1+n_1)\cdots Q(x_d+n_d)\\
&=&
\sum_{n_1=0}^{N-1} Q(x_1+n_1)
\cdots \sum_{n_d=0}^{N-1}Q(x_d+n_d)= 1.
\end{eqnarray*}
By construction, $P\chi_{[0,2]^d}\in C^{2L-1}(\mrd).$ This completes the proof.
\ep


\section{Construction of Gabor frames in $\ltrd$} \label{41226b}
As application of our results in Section \ref{75} and Section \ref{41226a}
we will now construct
pairs of dual Gabor frames in $\ltrd$ with attractive properties. Other constructions
in $\ltrd$
in the literature include \cite{CR} and \cite{KimII}; we will comment on these along the way.
Note also the approach (in $\ltr$) by Laugesen in \cite{Lau}.
In Section \ref{50401a} we will show that the constructions presented here yield
a very convenient way to obtain finite dual pairs of Gabor frames in $\ell^2(\mzd)$ as well.

For $\by\in \mrd$, the translation operator $T_{\by}$ and the
modulation operator $E_{\by},$ both acting on $\ltrd,$ are defined by
$(T_{\by}f)(\bx)  =  f(\bx-\by), \
(E_{\by}f)(\bx)  =  e^{2\pi i\by\cdot \bx}f(\bx), \ \ \bx\in
\mrd,$ where $\by\cdot \bx$ denotes the inner product of
$\by$ and $\bx$ in $\mrd$. Given a real and invertible $d\times
d$ matrix $B$ and $g\in \ltrd$ we consider Gabor systems of the form
$$\{E_{B\bm} T_\bn g\}_{\bm,\bn\in
\mzd} = \{ e^{2\pi iB \bm\cdot \bx} g(\bx-\bn)\}_{\bm,\bn\in
\mzd}.$$ The theory for duality of Gabor systems  is closely related to our discussion about partition of unities. Using the notation $B^\sharp = (B^T)^{-1},$ it is well known
(see \cite{RoSh1, Jan2, Lab,HLW}) that two Bessel sequences
$\{E_{B\bm}T_{\bn}g\}_{\bm,\bn\in \mzd}$ and
$\{E_{B\bm}T_{\bn}h\}_{\bm,\bn\in \mzd}$ form dual frames for
$\ltrd$ if and only if \bee \label{gfs} \sum_{\bk\in \mzd}
\overline{g(\bx- \bs \bn+\bk)}h(\bx+\bk)  = | \det
B|\delta_{\bn,0}, \, \mbox{a.e.} \, \bx\in \mrd,\ene for all $\bn \in \mzd.$
We will now choose the functions $g$ and $h$ of the form
\bes g(\bx)=G(\bx) \chi_{[0,N]^d}(\bx), h(\bx)=H(\bx) \chi_{[0,N]^d}(\bx), \bx\in \mrd,\ens
for some $N\in \mn$ and entire functions $G,H.$ Then, for $\bn=0,$ the condition \eqref{gfs}
takes the form
\bes \sum_{\bk\in \mzd}
\overline{G(\bx+\bk)}H(\bx+\bk)\chi_{[0,N]^d}(\bx+\bk)  = | \det B |, \, \, \mbox{a.e.} \, \bx\in \mrd;\ens up to the factor $| \det B |$ this is clearly a partition of unity condition on the
function $\overline{G}H\chi_{[0,N]^d}.$ Thus, the results in Section \ref{75} and Section \ref{41226a}
have almost immediate consequences for construction of dual Gabor frames.
For the discussion of the regularity of the obtained constructions, we will use that
if $P(\bx):= \prod_{j=1}^d \sin^{M}(\pi x_j/N)
A(\bx)$ for some real-valued and continuous $(N\mz)^d$-periodic
trigonometric polynomial $A$, then  $P\chi_{[0,N]^d} \in
C^{M-1}(\mr^d)$. The relation between the parameter $M$ in the following Theorem \ref{50105b}
and the parameter $L$ in Theorem \ref{eh-3} is that $L=2M.$

\bt \label{50105b} Let $K,N,M\in\mn$ with $M \le K.$ Assume that
\begin{equation*}
P(\bx)=\sum_{\bk \in \left(\mz \cap[-K,K]\right)^d} c_{\bk}
e^{2\pi i \bk \cdot \bx /N},\ \bx\in\mr^d
\end{equation*}
is a real-valued trigonometric polynomial with $
c_{\bk}=\frac{1}{N^d}\delta_{\bk,{\bf 0}}, \, \forall \bk\in
(N\mz)^d,$ and that $P\chi_{[0,N]^d}\in C^{2M-1}(\mr^d)$. Let $P$
be factorized as \bee \label{0803a} P(\bx)= \prod_{j=1}^d
\sin^{2M}(\pi x_j/N) G(\bx)H(\bx)\ene for some $(N\mz)^d$-periodic
real-valued trigonometric polynomials $G,H$. Let $B$  be a real
and invertible $d\times d$ matrix such that \bee \label{50105a}
B^{\sharp} \bn \notin (-N,N)^d , \ \forall \bn\in \mzd\setminus
\{{\bf 0}\}.\ene
 Then  the functions
\begin{eqnarray*}
g(\bx)&=& \left(\prod_{j=1}^d \sin^{M}(\pi x_j/N)\right)G(\bx)\chi_{[0,N]^d}(\bx), \label{eh-13}\\
h(\bx)&=& \left|  \det  B  \right| \left(\prod_{j=1}^d \sin^{M}(\pi
x_j/N)\right) H(\bx)\chi_{[0,N]^d}(\bx) \label{eh-14}
\end{eqnarray*}
belong to $ C^{M-1}(\mrd) $  and
 generate dual frames
$\{E_{B\bm}T_{\bn}  g\}_{\bm,\bn\in \mzd}$ and $\{E_{B\bm}T_{\bn} h\}_{\bm,\bn\in \mzd}$
for $\ltrd.$
\et
\bp
Since
$B^{\sharp} \bn \notin (-N,N)^d , \ \forall \bn\in \mzd\setminus
\{{\bf 0}\}$,
\eqref{gfs} is satisfied for $\bn\neq {\bf 0};$
thus the result follows from Theorem \ref{eh-3} and the comment just before the formulation
of Theorem \ref{50105b}.
\ep
Some comments are in order:
\bei
\item[(i)] Under the assumptions in Theorem \ref{50105b}, we know from Theorem \ref{eh-3}
that factorizations of $P$ as in \eqref{0803a} in terms of trigonometric polynomials
$G,H$ always exist. For example, we can choose $G(\bx)=\prod_{j=1}^d e^{2\pi i x_j M/N}$ and $H(\bx)=A_{2M}(\bx)$,
where $A_{2M}$ is defined as in \eqref{et-69}.

\item[(ii)] The construction in Theorem \ref{50105b} is much simpler than the one given in
\cite{CR}. We avoid a complicated book keeping; and we avoid to enlarge the support of the dual window and can keep the same support size as for the window itself. The construction by
I. Kim in \cite{KimI,KimII} also provides dual windows with the same support as the given window, but without an explicit expression for the dual window.
\item[(iii)] Small adjustments of  Theorem \ref{50105b} lead to constructions of tight frames.
If we (in addition to the stated conditions) assume that the trigonometric polynomial $P$ is non-negative,
the function
\begin{equation*}
 k(\bx):=\sqrt{\left| \det B  \right|}
\left(\prod_{j=1}^d \sin^{M}
(\pi x_j/N)\right)\sqrt{G(\bx)H(\bx)} \chi_{[0,N]^d}(\bx),
\end{equation*}
 belongs to $C^{M-1}(\mrd)$ and generates a tight Gabor frame $\{E_{B\bm}T_{\bn} k\}_{\bm,\bn\in \mzd}$ with frame bound $1.$
 \item[(iv)] A simple scaling extends Theorem \ref{50105b} to a construction of dual frames $\{E_{B\bm}T_{C\bn}  g\}_{\bm,\bn\in \mzd}$ and $\{E_{B\bm}T_{C\bn} h\}_{\bm,\bn\in \mzd},$ where $B$ and $C$ are real-valued invertible matrices. We leave the formulation of this
     result to the interested reader.
\eni

The following construction of a tight frame based on Corollary \ref{eh-5}
and the bullet (iii) above appears to be particularly simple and useful.

\bc \label{41228p}
Let $N\in \mn$ with $N\ge 2$, and
let $B$ be a real and invertible
$d\times d$ matrix such that \eqref{50105a} holds.
Let
$$ k(\bx) := \sqrt{\left| \det  B  \right|}\left( \frac{4^{N-1}}{ N {2N-2 \choose N-1} }\right)^{d/2}
\prod_{j=1}^d \sin^{N-1}(\pi x_j/N)
\chi_{[0,N]^d}(\bx) \label{eh-15}.$$
Then
$ k \in C^{N-2}(\mrd)$, and $\{E_{B\bm}T_{\bn} k\}_{\bm,\bn\in \mzd}$ is a tight Gabor frame with frame bound 1.
\ec

Note that this particular construction resembles the original work by
Daubechies, Grossmann, and Meyer \cite{DGM}: the difference is that
\cite{DGM} only deals with the case $N=2$ and the dimension $d=1.$ It is clear
that in Corollary
\ref{41228p} the possibility to increase $N$ is the key to the higher regularity.
We also note that  unless the matrix $B$ is a diagonal matrix, the construction
in Corollary
\ref{41228p} is not a tensor product of the constructions in \cite{DGM}.

Arbitrary (finite) regularity can be obtained for the windows $k$ in Corollary \ref{41228p}
by increasing the parameter $N\in \mn,$
but the price is that the support size increases as well.  The following alternative construction is based on the condition  \eqref{50105a} with  the choice $N=2,$
but nevertheless it allows us
to obtain arbitrary high regularity. The merits of this particular construction are twofold:
we are able to increase the regularity while we keep the support
$[0,2]^d,$ and the redundancy of the resulting Gabor frames is minimized
among the ones satisfying \eqref{50105a}  (the smallest values for $|\det \bs | = |\det B |^{-1}$ for matrices $B$ satisfying \eqref{50105a} are
clearly obtained whenever $N=2$.)
\bc\label{eh-15-4}
Let $L\in \mn$, and let $B$ be a real and invertible
$d\times d$ matrix such that \eqref{50105a} holds with $N=2.$
Define
$$ g(\bx)=  \prod_{j=1}^d \sin^{L}(\pi x_j/2)
\chi_{[0,2]^d}(\bx) $$
and
$$
h(\bx)=\left|  \det  B  \right| \prod_{j=1}^d \sin^{L}(\pi x_j/2)
\sum_{k_j=0}^{L-1}
{2L -1 \choose k_j}
   \sin^{2(L  -1 -k_j)}(\pi x_j/2) \cos^{2k_j} (\pi x_j/2) \chi_{[0,2]^d}(\bx).
$$
 Then $g \in
C^{L-1}(\mrd)$, $h \in C^{L-1}(\mrd)$, and the functions
$\{E_{B\bm}T_{\bn}  g\}_{\bm,\bn\in \mzd}$ and $\{E_{B\bm}T_{\bn} h\}_{\bm,\bn\in \mzd}$
 form
a pair of dual frames. \ec
\bp
\begin{figure}
\centerline{
\subfigure[]{\includegraphics[width=2.7in,height=2.5in]{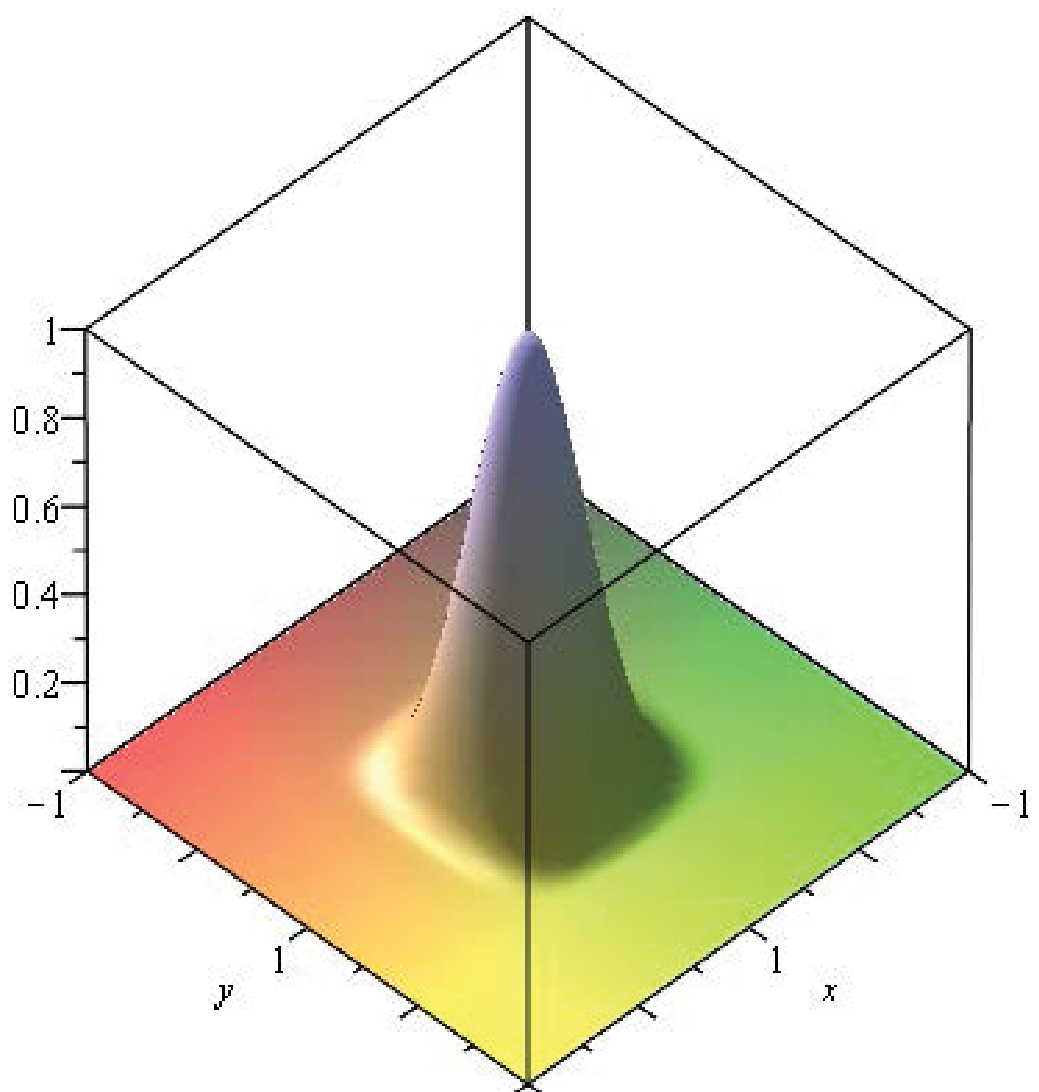}}\hfil
\subfigure[]{\includegraphics[width=2.7in,height=2.5in]{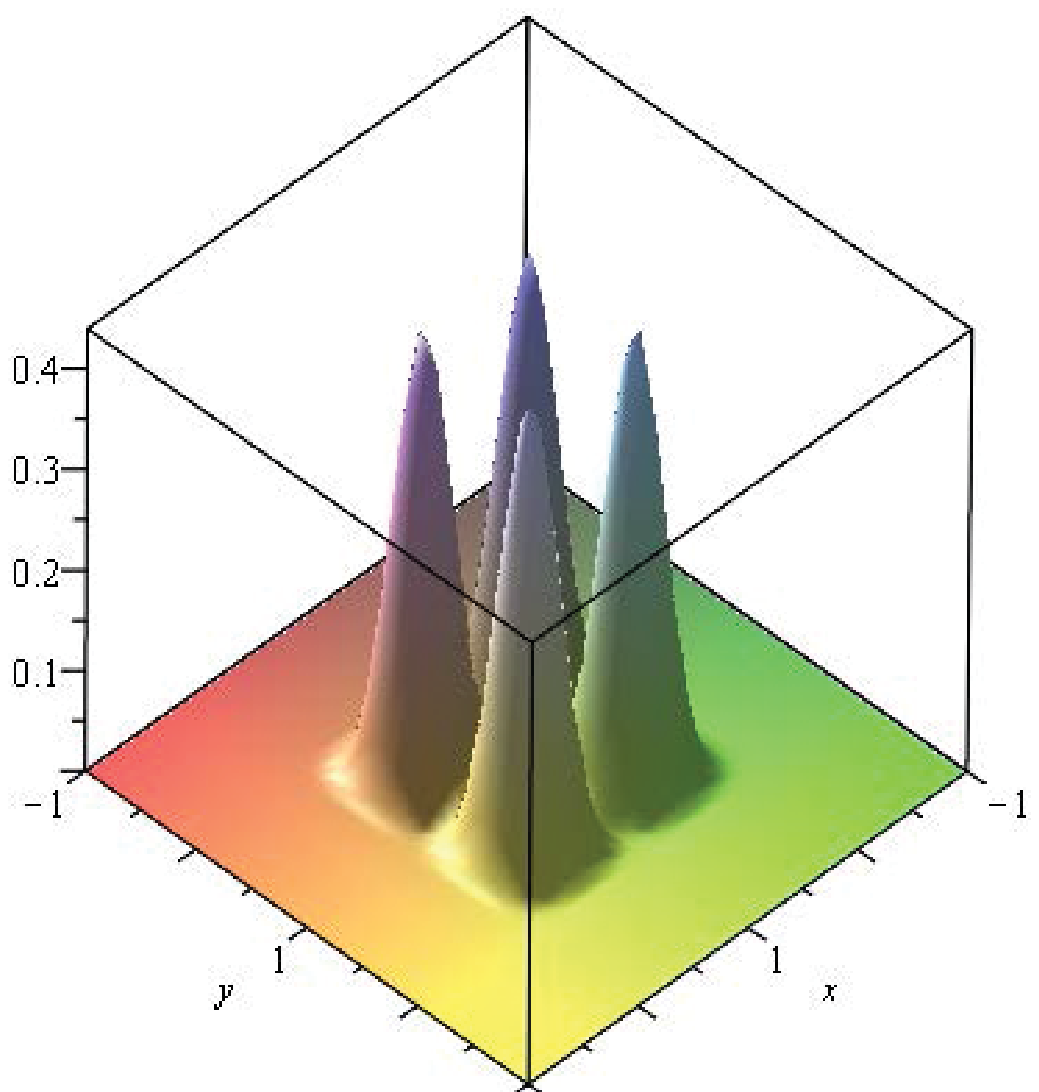}}}\hfil
\caption{The functions  $g$ and $h$ in Corollary \ref{eh-15-4} for
 $L=3$, $N=d=2$  and the matrix $B$ with entries $b_{11}=b_{22}=1/2, b_{12}=0, b_{21}=-1/2.$}
\label{eh-15-5} \end{figure} Using Corollary \ref{eh-15-3}, it
follows that the functions $g$ and $h$ satisfy the condition
\eqref{gfs} for $n=0.$ The choice of $B$ and the support sizes for
$g$ and $h$ shows that \eqref{gfs} holds for $n\neq 0$ as well.
\ep

We illustrate Corollary \ref{eh-15-4} in  Figure \ref{eh-15-5}, which is based on
the choices $L=3$, $N=d=2$  and $B=\frac{1}{2}\begin{pmatrix}
1 & 0\\
-1 & 1
\end{pmatrix}.$

We will now give a closer analysis of the central condition \eqref{50105a}.
For any $d\times d$ matrix, define the
norm $||B||$ by \bes ||B||= \sup_{||\bx||=1} \, ||B\bx||.\ens
We first state the following simple sufficient condition for
\eqref{50105a} to hold.

\bl\label{eh-10}  The condition \eqref{50105a} is satisfied if
$||B||\leq \frac{1}{\sqrt{d}\ N}$.
\el
\bp Since $B$ is
invertible, for any $\bn\in \mzd$ we have \bes ||\bn|| = || B^T
\bs \bn|| \le ||B|| \; ||\bs \bn||; \ens thus, for $\bn\neq 0$,
$||\bs \bn|| \ge ||\bn||/{||B||}\ge \sqrt{d} N$ if $||B||\leq
1/(\sqrt{d}N)$. Therefore, $ \bs \bn \notin (-N,N)^d$ for $\bn\in
\mzd\setminus \{{\bf 0}\} $. \ep

The following example shows that (at least for $N=2$) the norm condition $||B||\leq 1/(\sqrt{d}N)$
is optimal, in the sense that we for any number $a>1/(\sqrt{d}N)$ can find
an invertible matrix $B$ with $||B|| =a$ such that
\eqref{50105a} is not satisfied:

\begin{ex} \label{eh-11}
Let $N=d=2$. Then for any $\epsilon>0$, there exists an invertible $2\times 2$ matrix $B_\epsilon$ with
$||B_\epsilon||=  \frac{1+\epsilon}{2 \sqrt{2}}$ such that
$  \bs_\epsilon \bn \in (-2,2)^2$ for some $\bn\in \mz^2\setminus \{{\bf 0}\}$.
In fact,  consider $B_\epsilon=\dfrac{1+\epsilon}{4}\begin{pmatrix}
1 & 1\\
1 & -1
\end{pmatrix}.$
Then we have
\begin{eqnarray*}
|| B_\epsilon || &=& \sup_{\theta} \left|\left|\frac{1+\epsilon}{4}
\begin{pmatrix}
1 & 1\\
1 & -1
\end{pmatrix}
\begin{pmatrix}
\cos \theta \\
\sin \theta
\end{pmatrix}
\right|\right|  \\
&=& \sup_\theta  \frac{1+\epsilon}{4}
\sqrt{ (\cos\theta+\sin\theta)^2+(\cos\theta-\sin\theta)^2 }
=  \frac{1+\epsilon}{2\sqrt{2}}.
\end{eqnarray*}
But
$\bs_\epsilon (1,0)^T = \dfrac{2}{1+\epsilon}\begin{pmatrix}
1 & 1\\
1 & -1
\end{pmatrix}
\begin{pmatrix}
1 \\
0
\end{pmatrix}
=  \dfrac{2}{1+\epsilon} \begin{pmatrix}
1 \\
1
\end{pmatrix}
\in (-2,2)^2$. Thus \eqref{50105a} does not hold.
\ep

\end{ex}

Easy calculations show that in the special case of a real and invertible $d\times d$ diagonal matrix $B,$ the condition
\eqref{50105a} is satisfied if and only if $||B|| \le 1/N.$
On the other hand, the following example shows that \eqref{50105a} can be satisfied for non-diagonal matrices $B$ with arbitrary large
norm. Also in the work by I. Kim \cite{KimI} it was mentioned that
Gabor frame constructions with large matrix norm $||B||$ are possible.

\begin{ex} \label{eh-10}
Let $N=d=2$.
There exist invertible $2\times 2$ matrices
$B_a, a>0,$ such that $  \left(\bs_a \bn \right)\notin (-2,2)^2$ for all $ \bn\in
\mz^2\setminus \{{\bf 0}\}$ but $\lim_{a\rightarrow\infty}|| B_a ||
=\infty$. In
fact, consider $B_a=\dfrac{1}{4}\begin{pmatrix}
2 & 0\\
-2a & 2
\end{pmatrix}.$
Then
$\bs_a=\begin{pmatrix}
2 & 2a\\
0 & 2
\end{pmatrix}.$
Let $\bn=(n_1, n_2)^T\in \mz^2\setminus \{{\bf 0}\}$. Then
$\bs_a \bn=n_1 (2,0)^T +n_2(2a,2)^T.$ If $n_2 \neq 0$, then $\left(
\bs_a \bn \right)_2 \not\in (-2,2)$; if $n_2 = 0$ and $n_1 \neq 0$,
then $\left( \bs_a \bn \right)_1 \not\in (-2,2)$. Thus $  \left(\bs_a
\bn \right)\notin (-2,2)^2, \ \forall \bn\in \mz^2\setminus \{{\bf
0}\}$. Also,
\begin{eqnarray*}
|| B_a ||
&=& \sup_\theta  \frac{1}{4}
\sqrt{ (2\cos\theta)^2+(-2a\cos\theta+2\sin\theta)^2 } \ge \frac12 a.
\end{eqnarray*}
Hence $ \lim_{a \rightarrow \infty} || B_a || =\infty$. \ep
\end{ex}

Relating to the one-dimensional case, we know that if  $\{E_{mb}T_{n} g\}_{m,n\in\mz}$ is a Gabor frame $L^2(\mr)$,
then $0<b\le 1.$ The corresponding statement in higher dimensions is
that $|\det B| \le 1$  is a necessary
condition for $\{E_{B\bm}T_{\bn} g\}_{m,n\in\mz}$ to be a frame
for $\ltrd.$ In contrast to the one-dimensional case, the following example shows that \eqref{50105a} might not be satisfied for $d\ge 2,$ regardless how small $|\det B|$ is.

\bex Fix $\epsilon>0$ and consider for $a>0$ the matrix
$B_a=\epsilon^{1/2} \, \begin{pmatrix}
a & 0\\
0 & a^{-1}
\end{pmatrix}.$
Then $\det B_a= \epsilon$ for all $a>0,$ and
$||B_a||= \epsilon^{1/2} \max\{a,a^{-1}\}.$
Also,
$B_a^\sharp=\frac1{\epsilon^{1/2}} \, \begin{pmatrix}
a^{-1} & 0\\
0 & a
\end{pmatrix}.$
Clearly
$B_a^\sharp \begin{pmatrix} 1 \\ 0 \end{pmatrix} = \frac1{\epsilon^{1/2}}  \begin{pmatrix} a^{-1} \\ 0 \end{pmatrix}, \, \,  B_a^\sharp \begin{pmatrix} 0 \\ 1 \end{pmatrix} = \frac1{\epsilon^{1/2}}  \begin{pmatrix} 0 \\ a \end{pmatrix};$ thus, for an arbitrary value of $\epsilon$
we can find $a>0$ such that $\det B_a= \epsilon$ but for some $\bn \in \mz^2\setminus\{ {\bf 0} \},$
$B^\sharp \bn \in (-2, 2)^2.$ The example easily extends to any $\mrd, \, d\ge 2.$
\ep \enx

\section{Gabor frames in $\ell^2(\mz^d)$ through sampling} \label{50401a}
As further application of the results we will now show that
the attractive properties of the dual pairs of Gabor frames in Section \ref{41226b}
yield an  easy way to construct Gabor frames in $\ell^2(\mz^d)$ through sampling. For general functions in $\ltr$ the work by Janssen \cite{Jan6} shows that sampling is a delicate
issue, but the continuity and compact support of the windows constructed in
Section \ref{41226b} remove several technical difficulties.
For further
information on discrete Gabor systems we refer to the paper \cite{Jan2} by Janssen (which
also deals with the more general case of shift-invariant systems),
\cite{CV1,CV2} by Cvetkovi\'c and Vetterli, as well as to the recent paper \cite{LoHa} by
Lopez and Han. We also mention that the theory for translation invariant systems on LCA groups
yields a joint framework to Gabor theory on $\ltrd$ and $\ell^2(\mz^d),$ see, e.g., the papers
\cite{LeSi1, LeSi2} by Jakobsen and Lemvig.

Let $B$ denote an invertible $d\times d$ matrix for which $B^{-1}$ has integer entries.
Consider the subgroup $G:= B^{-1} \mzd$ of $\mzd,$ and let $\Omega$ denote a collection of
coset representatives of the coset $\mzd/G;$ that is, $\mzd$ is a disjoint union of the
sets $G+ \bm,$ where $\bm \in \Omega.$ It is well known that the number of elements in
$\Omega$ is
\bes | \Omega| = | \det (B^{-1}) | = \frac1{| \det B |}.\ens
Let us now fix any sequence in $\ell^2(\mzd);$ for our current purpose it will be convenient to denote the
sequence by $\{c(\bj)\}_{\bj\in \mzd}$ rather than $\{c_{\bj}\}_{\bj\in \mzd}.$ We will consider the {\it Gabor system in $\ell^2(\mzd)$ generated by the sequence $\{c(\bj)\}_{\bj\in \mzd}$ and
the matrix $B,$} i.e., the collection of sequences $\{c_{\bm,\bn}\}_{\bm,\bn\in \mzd}\subset \ell^2(\mzd)$
given by
\begin{equation*} 
c_{\bm,\bn}(\bj)= e^{2\pi i B\bm \cdot \bj} c(\bj-\bn), \, \bj\in \mzd.
\end{equation*}
Given two sequences $\{c(\bj)\}_{\bj\in \mz}, \{d(\bj)\}_{\bj\in \mz} \in \ell^2(\mzd)$
such that $\{c_{\bm,\bn}\}_{\bm,\bn\in \mzd}$ and $\{d_{\bm,\bn}\}_{\bm,\bn\in \mzd}$
are Bessel sequences, it was shown in Theorem 1.4 in \cite{LoHa} that
$\{c_{\bm,\bn}\}_{\bm,\bn\in \mzd}$ and $\{d_{\bm,\bn}\}_{\bm,\bn\in \mzd}$ are dual frames for
$\ell^2(\mzd)$ if and only if
\bee \label{50330a} \sum_{k\in \mzd} \overline{c(\bj - \bs \bn + \bk)} d(\bj + \bk)= | \det B| \, \delta_{\bn,0}\ene for all $\bj, \bn \in \mzd.$

The relation to the results in Section \ref{41226b} is evident. In fact, for all the
frame constructions in Section \ref{41226b}, the windows $g(x)=G(x) \chi_{[0,N]}(x)$ and
the dual windows $h(x)=H(x) \chi_{[0,N]}(x)$ are continuous functions with compact support, which implies that the duality conditions \eqref{gfs} hold {\it pointwise} for all $x\in \mrd;$
this clearly implies that the sequences $\{c(\bj)\}_{\bj\in \mzd}:= \{g(\bj)\}_{\bj\in \mzd}$ and $\{d(\bj)\}_{\bj\in \mzd}:=\{h(\bj)\}_{\bj\in \mzd}$ satisfy \eqref{50330a}.
In other words: the samples of the dual windows $g,h,$ i.e., the sequences $\{g(\bj)\}_{\bj\in \mzd}$ and $
\{h(\bj)\}_{\bj\in \mzd}$ generate dual Gabor frames for $\ell^2(\mzd).$
Thus, Theorem \ref{50105b} has the following immediate consequence:

\bc Under the assumptions in Theorem \ref{50105b}, the sequences
\begin{eqnarray*}
g({\bf j})&:=& \left(\prod_{\ell=1}^d \sin^{M}(\pi {\bf j}_\ell/N)\right)G({\bf j})\chi_{[0,N]^d}({\bf j}), \, \bj \in \mzd, \\
h({\bf j})&=& \left|  \det  B  \right| \left(\prod_{\ell=1}^d
\sin^{M}(\pi {\bf j}_\ell/N)\right) H({\bf j})\chi_{[0,N]^d}({\bf j}), \, \bj \in \mzd,
\end{eqnarray*}
belong to $\ell^2(\mzd)$  and the associated discrete Gabor systems
$\{g_{\bm,\bn}\}_{\bm,\bn\in \mzd}$ and
$\{h_{\bm,\bn}\}_{\bm,\bn\in \mzd}$ form dual frames for  $\ell^2(\mzd)$. \ec

Similarly, Corollary \ref{41228p} leads to the following explicit constructions of
finite tight Gabor frames in $\ell^2(\mzd):$

\bc \label{eh-32} Let $N\in \mn$ with $N\ge 2$, and let $B$ be a
real and invertible $d\times d$ matrix for which $B^{-1}$ has
integer entries and \eqref{50105a} holds. Let
$$ k({\bf j}) :=
\sqrt{\left| \det  B  \right|}\left( \frac{4^{N-1}}{ N {2N-2
\choose N-1} }\right)^{d/2} \prod_{\ell=1}^d \sin^{N-1}(\pi {\bf
j}_\ell/N) \chi_{[0,N]^d}({\bf j}), \, \bj\in \mzd.$$
Then $
\{k({\bf j})\}_{\bj \in \mzd} \in \ell^2(\mzd)$, and the
associated discrete Gabor system $\{k_{\bm,\bn}\}_{\bm,\bn\in
\mzd}$
 is a tight Gabor frame for $\ell^2(\mzd)$ with frame bound 1.
\ec

Note that sampling of the frame in Corollary \ref{eh-15-4} yields the trivial
dual pair of Gabor frames in $\ell^2(\mzd),$ with windows vanishing at all other points than
$(1,1, \dots,1).$

\noindent{\bf Acknowledgment:} The first-named author would like to thank Peter Massopust for
useful discussions about the content of the manuscript. The authors also thank Hans Feichtinger for suggesting to include the application to discrete Gabor frames. They also thank the reviewers for
their detailed comments, which improved the manuscript.

{\bf \vspace{.1in}

\noindent Ole Christensen\\
Department of Applied Mathematics and Computer Science\\
Technical University of Denmark, Building 303,
2800 Lyngby  \\
Denmark \\
 Email: ochr@dtu.dk

\vspace{.1in}\noindent Hong Oh Kim \\
Division of General Studies, UNIST\\
50 UNIST-gil, Ulsan 44919\\
Republic of Korea\\
Email: hkim2031@unist.ac.kr

\vspace{.1in} \noindent Rae Young Kim \\
Department of Mathematics,
Yeungnam University\\
280 Daehak-Ro, Gyeongsan, Gyeongbuk 38541\\
Republic of Korea\\
Email:  rykim@ynu.ac.kr }

\end{document}